\documentclass[reqno,12pt]{amsart}
\usepackage{amsmath,amsthm,amssymb,amsfonts,amscd}
\usepackage{graphicx,color}
\usepackage{boxedminipage}
\usepackage{array}

\input xy
\xyoption{all}

\newcommand{\Z}{{\mathbb{Z}}}

\def \cA{{\mathcal A}}

\def \cC{{\mathcal C}}
\def \cD{{\mathcal D}}
\def \cE{{\mathcal E}}

\def \cX{\mathcal{X}}
\def \cY{\mathcal{Y}}

\def \Id{\mathrm{Id}}

\def \Im{\mathrm{Im}}
\def \Ker{\mathrm{Ker}}

\def \Ob{\mathrm{Ob}}

\def \Tw{\mathit{Tw}}
\def \Tri{\mathit{Tr}}

\def \1{{\bf 1}}
\def \b1{{\bar {\bf 1}}}
\def \b{{\frak b}}
\def \f{{\frak f}}
\def \g{{\frak g}}

\def \m{{\frak m}}

\def \ti{\tilde}

\newcommand{\bp}{\begin{pmatrix}}
\newcommand{\ep}{\end{pmatrix}}

\newcommand{\bps}{\begin{smallmatrix}}
\newcommand{\eps}{\end{smallmatrix}}

\def \({\left(}
\def \){\right)}
\def \tri{\triangle}

\def \0{{\bf 0}}
\def \1{{\bf 1}}

\def \f{{\frak f}}

\def \Nsddata#1#2#3#4#5{
(\xymatrix{{#1}\  \ar@<0.5ex>[r]^{{#2}} & \ {#4}
\ar@<0.5ex>[l]^{{#3}}} ,#5) }

\def \Nsddata#1#2#3#4#5{
(\xymatrix{{#1}\  \ar@<0.5ex>[r]^{{#2}} & \ {#4}
\ar@<0.5ex>[l]^{{#3}}} ,#5) }

\def \HD#1#2#3#4{
\xymatrix{{#1}\  \ar@<0.5ex>[r]^{{#2}} & \ {#4}
\ar@<0.5ex>[l]^{{#3}}}}

\newtheorem{thm}{Theorem}[section]

\newtheorem{lem}[thm]{Lemma}
\newtheorem{prop}[thm]{Proposition}

\theoremstyle{definition}

\newtheorem{defn}[thm]{Definition}%[section]

\begin{document}

\title{Comments on the shifts and the signs in $A_\infty$-categories} 

%\date{\today}

\noindent{
 }\\

\author{Hiroshige Kajiura}
\address{Faculty of Science, Chiba university, 263-8522, Japan}
\email{kajiura@math.s.chiba-u.ac.jp}

\begin{abstract}
We discuss that there exist at least
two different choices in the 
signs of the induced $A_\infty$-structures 
in shifting the degree of objects in an $A_\infty$-category. 
We show that both of these choices are natural
in the sense that they are compatible
with the homological perturbation theory
of $A_\infty$-categories.
These different choices give different triangulated categories 
out of an $A_\infty$-category. 
\end{abstract}

\renewcommand{\thefootnote}{}
\footnote{This work 
is supported by Grant-in Aid for Scientific Research (C) (18K03293) 
%(25400081)
of the Japan Society for the Promotion of Science. }

\renewcommand{\thefootnote}{\arabic{footnote}}
\setcounter{footnote}{0}

\maketitle

{\small



}

\tableofcontents

%%%%%%%%%%%%%%%%%%%%%%%%%%%%%%%%%%%%%%%%%%%%%%%%%%%%%%

%\hspace*{7.0cm} \rotatebox[origin=c]{90}{\green{$=$}}
%\rotatebox[origin=c]{20}{\scalebox{1}[2]{A}}

 \section{Introduction}

In order to formulate his homological mirror symmetry conjecture, 
Kontsevich proposes a way of constructing a triangulated category $\Tri(\cC)$ 
from an $A_\infty$-category $\cC$ \cite{HMS}. 
This is a natural generalization of 
the Bondal-Kapranov's DG-enhancement \cite{BoKa:enhanced} 
to $A_\infty$ setting, and 
consists of three steps. 
First, for a given $A_\infty$-category $\cC$, 
we add formally shifted objects of $\cC$ and their direct sums
and construct 
an additive enlargement $\ti{\cC}$.
We then construct a further enlarged $A_\infty$-category $\Tw(\cC)$ 
consisting of one-sided twisted complexes, 
and the zero-th cohomology $\Tri(\cC):=H^0(\Tw(\cC))$ 
turns to form a triangulated category. 
\footnote{Precisely speaking, 
what we call $A_\infty$-categories here should be unital 
ones in some sense as in \cite{lyuba:unital}. }

In this paper, we discuss that there exist at least two different
natural choices 
of inducing $A_\infty$-structures on $\ti{\cC}$ from that in $\cC$. 
In order to simplify the signs appearing the defining equations of
$A_\infty$-structures, we usually express them
in the suspended notation.
In this suspended notation, 
a natural induced $A_\infty$-structure in $\ti{\cC}$ from $\cC$
is given in \cite{fukaya:mirror2}. 
Let $\{m_k\}_{k\ge 1}$ be the $A_\infty$-structure in $\cC$, and 
we denote by $\{b_k:=s m_k (s^{-1})^{\otimes k}\}_{k\ge 1}$ the one
in the suspended notation. 
For $X_i,X_j\in\cC$ and $\alpha_{ij}\in\cC(X_i,X_j)$, 
we denote $\bar{\alpha}_{ij}:=s(\alpha_{ij})\in s\cC(X_i,X_j)$,
and denote by $\bar{\alpha}'_{ij}\in s\cC(X_i[r_i],X_j[r_j])$
the corresponding morphism. 
Then, one can easily check that $\{\ti{b}_k\}$ defined by 
\begin{equation*}
  \ti{b}_k(\bar{\alpha}'_{12}, \dots, \bar{\alpha}'_{k(k+1)})
  := (-1)^{r_1}\, (b_k(\bar{\alpha}_{12},\dots,\bar{\alpha}_{k(k+1)}) )' , 
\end{equation*}
where $\alpha_{i(i+1)}\in\cC(X_i,X_{i+1})$, 
preserves the defining relations of $A_\infty$-structure. 

However, there are other natural choices of $\{\ti{b}_k\}$,
or equivalently, $\{\ti{m}_k\}$. 
One such is defined by 
\begin{equation*}
  \ti{m}_k(\alpha'_{12}, \dots, \alpha'_{k(k+1)})
  := (-1)^{r_1k}\, (m_k(\alpha_{12},\dots,\alpha_{k(k+1)}) )' . 
\end{equation*}
In the suspended notation, this turns out to be
\begin{equation*}
  \ti{b}_k(\bar{\alpha}'_{12}, \dots, \bar{\alpha}'_{k(k+1)})
  := (-1)^{r_1+\cdots + r_k}\,
  (b_k(\bar{\alpha}_{12},\dots,\bar{\alpha}_{k(k+1)}) )' . 
\end{equation*}
We denote the former $A_\infty$-category and the latter $A_\infty$-category
by $\ti{\cC}^{(1)}$ and $\ti{\cC}^{(2)}$, respectively. 
In general, these two cannot be linear $A_\infty$-isomorphic to each other. 
We discuss that both of these two 
are natural by showing
that both of these preserve $A_\infty$-morphisms in the sense in 
Proposition \ref{prop:morp-preserve} 
and preserve the homological perturbation theory (HPT) 
in the sense of Theorem \ref{thm:HPT-preserve}. 
It is a typical property of homotopy algebras such as $A_\infty$-algebras
that the homological perturbation theory works well
as developed in \cite{gugjds,GLS:chen,GLS,HK}, etc.
In this sense, we believe that Theorem \ref{thm:HPT-preserve}
guarantees
that both of these two choices are natural enough. 

Let us explain how the latter choice $\ti{\cC}^{(2)}$ also
arises naturally. 
Recall that an $A_\infty$-category $\cC$ with higher products all zero,
$m_3=m_4=\cdots=0$, is a DG category. 
In particular,
for a given additive category $\cA$, 
the DG category $DG(\cA)$ of complexes in $\cA$ has its differential
and its composition of the type $\ti{\cC}^{(2)}$. 
\footnote{More generally, we may begin with
a {\em full subcategory $\cA$ of} an additive category and then obtain 
the DG category $DG(\cA)$. }
Namely, 
instead of considering the category $Comp(\cA)$ of complexes in $\cA$, 
we can naturally construct a DG category $DG(\cA)$
whose objects are the same as those in $Comp(\cA)$. 
For two complexes $X^\bullet$, $Y^\bullet$ in $\cA$, 
a map in 
\begin{equation*}
 DG^r(\cA)(X^\bullet,Y^\bullet)  
\end{equation*}
is a collection $\varphi=\{\varphi^i:X^i\to Y^{i+r}\}_{i\in\Z}$
of homomorphisms. 
Thus, $DG(\cA)$ has a natural associative composition of morphisms,
which we denote by $\varphi\circ\psi$
for $\varphi\in DG^r(\cA)(X^\bullet,Y^\bullet)$ and
$\psi\in DG^r(\cA)(Y^\bullet,Z^\bullet)$. 
The differential 
\begin{equation*}
 d: DG^r(\cA)(X^\bullet,Y^\bullet)\to DG^{r+1}(\cA)(X^\bullet,Y^\bullet) 
\end{equation*}
is then defined by
\begin{equation*}
 d(\varphi) := d_X \circ \varphi - (-1)^r \varphi\circ d_Y, 
\end{equation*}
or more explicitly, by 
\begin{equation*}
( d(\varphi) )^{i} := d_X^i \circ \varphi^{i+1} - (-1)^r\varphi^i\circ d_Y^{i+r} 
\end{equation*}
for $\varphi\in DG^r(\cA)(X^\bullet,Y^\bullet)$. 
When we denote
\begin{equation*}
 \begin{split}
   Z^r(DG(\cA))(X^\bullet,Y^\bullet)& :=\Ker (d: DG^r(\cA)(X^\bullet,Y^\bullet)\to DG^{r+1}(\cA)(X^\bullet,Y^\bullet) ), \\
   B^r(DG(\cA))(X^\bullet,Y^\bullet)& :=\Im (d: DG^{r-1}(\cA)(X^\bullet,Y^\bullet)\to DG^{r}(\cA)(X^\bullet,Y^\bullet) ),   
 \end{split}
\end{equation*}
the category $Z^0(DG(\cA))$ is nothing but $Comp(\cA)$, 
whose homotopy category is then the zero-th cohomology
$H^0(DG(\cA))$. 

Here, the shift functor $T:DG(\cA)\to DG(\cA)$ is defined by
$T(X^\bullet)=(X[1])^\bullet=(\{(X[1]^i,d_{X[1]}^i\}_{i\in\Z})$, where
\begin{equation*}
  X[1]^i= X^{i+1} ,\qquad d_{X[1]}^i = -d_X^{i+1} . 
\end{equation*}
Then, we can identify an element $\varphi\in DG^r(\cA)(X^\bullet,Y^\bullet)$
with $\mu\in DG^{r-1}(\cA)(X^\bullet,Y[1]^\bullet)$ by
\begin{equation*}
  \mu^{i}=\varphi^i . 
\end{equation*}
Similarly, we can identify an element $\varphi\in DG^r(\cA)(X^\bullet,Y^\bullet)$
with $\nu\in DG^{r+1}(\cA)(X[1]^\bullet,Y^\bullet)$ by
\begin{equation*}
  \nu^{i-1}=\varphi^i . 
\end{equation*}

Let us discuss compatibility of these identifications with
the structure of the DG category. 
For $\mu\in DG^{r-1}(\cA)(X^\bullet,Y[1]^\bullet)$, 
the differential acts as
\begin{equation*}
  \begin{split}
    ( d(\mu) )^{i}
    & = d_X^i \circ \mu^{i+1} - (-1)^{r-1}\mu^i\circ d_{Y[1]}^{i+r-1} \\
    & = d_X^i \circ \varphi^{i+1} - (-1)^r\varphi^i\circ d_{Y}^{i+r} \\
    & = (d(\varphi))^i .
  \end{split}
\end{equation*}
On the other hand, 
for $\nu\in DG^{r+1}(\cA)(X[1]^\bullet,Y^\bullet)$, 
the differential acts as
\begin{equation*}
  \begin{split}
    ( d(\nu) )^{i-1}
    & = d_{X[1]}^{i-1} \circ \nu^{i}
    - (-1)^{r+1}\nu^{i-1}\circ d_{Y}^{(i-1)+(r+1)} \\
    & = -d_{X}^{i} \circ \varphi^{i+1} + (-1)^{r}\varphi^i\circ d_{Y}^{i+r} \\
    & = -(d(\varphi))^i .
  \end{split}
\end{equation*}
The composition of morphisms is compatible with the identifications
in the sense that we have no sign. 

Let us denote
$d=:\ti{m}_1$ and $\varphi\circ\psi=:\ti{m}_2(\varphi,\psi)$ 
for $\varphi\in DG(\cA)(X^\bullet, Y^\bullet)$ and
$\psi\in DG(\cA)(Y^\bullet, Z^\bullet)$. 
Furthermore, for $r_1,r_2,r_3\in\Z$, we denote by 
$\varphi'\in DG(\cA)(X[r_1]^\bullet, Y[r_2]^\bullet)$ and
$\psi'\in DG(\cA)(Y[r_2]^\bullet, Z[r_3]^\bullet)$
the morphism identified with $\varphi$ and $\psi$, respectively. 
Then, the compatibilities are summarized as
\begin{equation*}
  \ti{m}_1(\varphi')=(-1)^{r_1} (\ti{m}_1(\varphi))',\qquad
  \ti{m}_2(\varphi',\psi')= (\ti{m}_2(\varphi,\psi))' . 
\end{equation*}
This implies that $DG(\cA)=\ti{\cC}^{(2)}$ by choosing 
a full subcategory $\cC\subset DG(\cA)$ suitably. 

We cannot know from this example
how we should determine the compatibilities for higher products. 
However, applying the HPT to a DG category such as $DG(\cA)$, 
it turns out that all the derived higher products
are exactly of the $\ti{\cC}^{(2)}$ type. 
This fact corresponds to the case that $\cC$ is a DG-category
in Theorem \ref{thm:HPT-preserve}. 
In this sense, we may say
that the latter choice $\ti{\cC}^{(2)}$ is natural in the 
DG notation or the unsuspended notation. 

In the next section, we recall terminologies
of $A_\infty$-algebras and $A_\infty$-categories 
in order to fix the notations, and also recall the HPT. 
Then, in section \ref{sec:main} we define $\ti{\cC}^{(a)}$, $a=1,2$,
and discuss these properties. 
In particular, we show 
Proposition \ref{prop:morp-preserve} and Theorem \ref{thm:HPT-preserve} hold. 
Once $\ti{\cC}^{(a)}$, $a=1$ or $2$, is constructed, 
one can further construct the $A_\infty$-category $\Tw^{(a)}(\cC)$ of
one-sided twisted complexes, and the triangulated category
$\Tri^{(a)}(\cC):=H^0(\Tw^{(a)}(\cC))$. 
These two triangulated categories $\Tri^{(1)}(\cC)$ and $\Tri^{(2)}(\cC)$
are not isomorphic to each other in any natural way. 
In subsection \ref{ssec:tri}, 
we present the shift functor $T:\Tri^{(a)}(\cC)\to\Tri^{(a)}(\cC)$ 
since it may not be presented explicitly in literatures.

%\vspace*{0.5cm}
%
%{\bf Acknowledgements}\quad 

 \section{$A_\infty$-categories and their properties}
\label{sec:recall}

In order to fix notations, 
we recall basic terminologies 
and properties of $A_\infty$-algebras 
and $A_\infty$-categories. 
Throughout this paper, 
all vector spaces are those over a fixed base field $K$.

\subsection{$A_\infty$-algebras}
\label{ssec:Ainfty}

\begin{defn}[$A_\infty$-algebra \cite{jds:hahI,jds:hahII}]
\label{defn:Ainfty}
An {\em $A_\infty$-algebra} $(A,\m)$ consists of a 
$\Z$-graded vector space $A:=\oplus_{r\in\Z} A^r$ 
with a collection of 
multilinear maps $\m:=\{m_n:A^{\otimes n}\to A\}_{n\ge 1}$ of 
degree $(2-n)$ satisfying 
\begin{equation}\label{Ainfty-rel}
  0=\sum_{k+l=n+1}\sum_{j=0}^{k-1}
(-1)^\star  
 \ m_k(a_1,\dots,a_j, m_l(a_{j+1},\dots,a_{j+l})
 ,a_{j+l+1},\dots,a_n) 
\end{equation}
for $n\ge 1$ with homogeneous elements $a_i\in A^{|a_i|}$, $i=1,\dots,n$, 
where the sign is given by
$\star=(j+1)(l+1)+l(|a_1|+\cdots+|a_j|)$. 
\end{defn}
That the multilinear map $m_k$ has degree $(2-k)$ indicates 
the degree of $m_k(a_1,\dots,a_k)$ is $|a_1|+\cdots +|a_k|+(2-k)$. 
The $A_\infty$-relations (\ref{Ainfty-rel}) imply $(m_1)^2=0$ for $n=1$, 
the Leibniz rule of the differential $m_1$ with respect to the product $m_2$ 
for $n=2$, and the associativity of $m_2$ up to homotopy for $n=3$. 
These facts further imply that 
the cohomology $H(A):=H(A,m_1)$ has the structure of 
a (non-unital) graded algebra, where 
the product is induced from $m_2$. 
We denote this algebra by $H(A,\m)$ and call 
the {\em cohomology algebra} of $(A,\m)$.

Note that 
the product $m_2$ is strictly associative in $A$ if $m_3$=0. 
\begin{defn}\label{defn:DGA}
An $A_\infty$-algebra $(A,\m)$ with higher products all zero, 
$m_3=m_4=\cdots =0$, is called 
a {\em differential graded (DG) algebra}. 
\end{defn}

Let $s: A^r\to (A[1])^{r-1}$ be the suspension. 
We denote the induced $A_\infty$-structure by
$\b=\{b_n:(A[1])^{\otimes k}\to A[1]\}_{n\ge 1}$. 
Explicitly, $b_k$ is defined by 
\begin{equation}\label{m-suspension}
 \begin{split}
 & b_k (s(a_1),\dots,s(a_k)) 
 =(-1)^\star s\cdot m_k(a_1,\dots,a_k), \\
 & \star = {(k-1)|s(a_1)|+(k-2)|s(a_2)|+ \cdots +2|s(a_{k-2})|
+ |s(a_{k-1})|} , 
 \end{split}
\end{equation}
where the degree of $b_k$ is one for any $k$. 
Hereafter, we denote 
$\bar{a}:=s(a)\in A[1]$. 
For the $A_\infty$-algebra $(A[1],\b)$ in this suspended notation, 
the sign in the $A_\infty$-relation (\ref{Ainfty-rel}) 
is simplified as 
\begin{equation}\label{Ainfty-rel2}
 \begin{split}
 & 0= \sum_{k+l=n+1}\sum_{j=0}^{k-1}
(-1)^\star
 \ b_k(\bar{a}_1,\dots,\bar{a}_j, b_l(\bar{a}_{j+1},\dots,\bar{a}_{j+l})
 ,\bar{a}_{j+l+1},\dots,\bar{a}_n) \\
 & \star= {|\bar{a}_1|+\cdots +|\bar{a}_j|} 
 \end{split}
\end{equation}
(Getzler-Jones \cite{gj:cyclic}). 

Let $T^c(A[1]):=\oplus_{k\ge 0} (A[1])^{\otimes k}$ 
be the tensor coalgebra of $A[1]$, where $(A[1])^{\otimes 0}:=K$. 
Here, the coproduct 
$\tri: T^c(A[1])\to T^c(A[1])\bigotimes T^c(A[1])$ 
is defined by 
\begin{equation*}
 \tri(\bar{a}_1\otimes\cdots\otimes \bar{a}_n)
 :=\oplus_{j=0}^n (\bar{a}_1\otimes\cdots\otimes \bar{a}_j)
 \bigotimes (\bar{a}_{j+1}\otimes\cdots\otimes \bar{a}_n)
\end{equation*}
for $\bar{a}_1\otimes\cdots\otimes \bar{a}_n\in (A[1])^{\otimes n}$, 
where, for instance, the term of $j=0$
reads $\1_K\bigotimes (\bar{a}_1\otimes\cdots\otimes \bar{a}_n)$. 
The $A_\infty$-structure $b_k$ 
is lifted to a coderivation $\b_k:T^c(A[1])\to T^c(A[1])$, 
$\tri\circ\b_k=(\b_k\otimes id_{T^c(A[1])}+ id_{T^c(A[1])}\otimes\b_k)\tri$. 
For the degree one 
coderivation $\b:=\b_1+\b_2+\cdots$, 
$(\b)^2=0$ follows from the $A_\infty$-relation (\ref{Ainfty-rel2}). 
Thus, 
an $A_\infty$-algebra $(A,\b)$ is equivalent to 
a DG coalgebra $(T^c(A[1]),\b,\tri)$. 
This construction is called 
the bar construction of $(A,\m)$. 

For two $\Z$-graded vector spaces $A,B$, 
coalgebra maps $\f: T^c(A[1])\to T^c(B[1])$, 
$\tri\f=(\f\bigotimes \f)\tri$, $\f(\1_K)=\1_K$, are in one-to-one 
correspondence with collections of degree zero multilinear maps 
$\{f_i: (A[1])^{\otimes i}\to B[1]\}_{i=1,\dots}$. 
In fact, for 
$\bar{a}_1\otimes\cdots\otimes \bar{a}_n\in (A[1])^{\otimes n}$, 
$\f$ is defined by 
\begin{equation*}
 \begin{split}
 & \f(\bar{a}_1\otimes\cdots\otimes \bar{a}_n) 
 :=\sum_{i\ge 1} \sum_{j_1+\cdots +j_i=n} \\
 &\quad f_{j_1}(\bar{a}_1,\dots,\bar{a}_{j_1})\otimes 
 f_{j_2}(\bar{a}_{j_1+1},\dots, \bar{a}_{j_1+j_2})\otimes\cdots 
 \otimes 
  f_{j_i}(\dots,\bar{a}_n) . 
 \end{split}
\end{equation*}
\begin{defn}[$A_\infty$-morphism]
For two $A_\infty$-algebras $(A,\m^A)$, $(B,\m^B)$, 
an {\em $A_\infty$-morphism} $\f:(A,\m^A)\to (A',\m^B)$ is 
a collection $\{f_i: (A[1])^{\otimes i}\to B[1]\}_{i=1,2,\dots}$ 
of degree zero multilinear maps 
whose lift $\f:T^c(A[1])\to T^c(B[1])$ to a coalgebra map 
gives a morphism $\f:(T^c(A[1]),\b^A)\to T^c(B[1],\b^B)$ 
between the complexes:  
\begin{equation*}
 \b^B\circ\f =\f\circ\b^A . 
\end{equation*}
In particular, we call an $A_\infty$-morphism $\f$ with 
$f_2=f_3=\cdots =0$ a {\em linear} $A_\infty$-morphism. 
\end{defn}
The condition $\b^B\circ\f =\f\circ\b^A$ is equivalent to that 
$\f$ satisfies the following equations 
\begin{equation}\label{Ainfty-morp}
 \begin{split}
 & \sum_{i\ge 1} \sum_{j_1+\cdots +j_i=n}
   b^B_i(f_{j_1}(\bar{a}_1,. .,\bar{a}_{j_1}),
   f_{j_2}(\bar{a}_{j_1+1},. ., \bar{a}_{j_1+j_2}),\dots,
   f_{j_i}(\dots,\bar{a}_n)) \\
 & = \sum_{k+l=n+1} \sum_{j=0}^{k-1} (-1)^{|\bar{a}_1|+\cdots + |\bar{a}_j|} 
   f_k(\bar{a}_1, \dots,\bar{a}_j,
   b^A_l(\bar{a}_{j+1},. ., \bar{a}_{j+l})\dots, \bar{a}_n) 
 \end{split}
\end{equation}
for any $n\ge 1$ and homogeneous elements $\bar{a}_1,\dots, \bar{a}_n\in A[1]$. 
This defining equation (\ref{Ainfty-morp}) of 
the $A_\infty$-morphism for $n=1$ implies that 
$f_1:A[1]\to B[1]$ forms a chain map $f_1:(A[1],m^A_1)\to (B[1],m^B_1)$. 
This together with 
the defining equation (\ref{Ainfty-morp}) for $n=2$ then implies that 
$f_1:A[1]\to B[1]$ induces a (non-unital) graded algebra map 
from $H(A,\m^A)$ to $H(B,\m^B)$. 
%We denote it by $H(\f):H(A,\m^A)\to H(B,\m^B)$. 
%

The composition of two $A_\infty$-morphisms is defined by 
the composition of the corresponding two coalgebra maps. 
Thus, the composition is associative. 
\begin{defn}
An $A_\infty$-morphism $\f:(A,\m^A)\to (B,\m^B)$ is 
called an {\em $A_\infty$-quasi-isomorphism}  
iff $f_1:(A[1],b^A_1)\to (B[1],b^B_1)$ 
induces an isomorphism between 
the cohomologies of these two complexes. 
If in particular $f_1$ is an isomorphism of complexes, 
then $\f$ is called an {\em $A_\infty$-isomorphism}. 
\end{defn}

 \subsection{$A_\infty$-categories}
\label{ssec:Ainftycat}

\begin{defn}[$A_\infty$-category \cite{fukayaAinfty}]
An {\em $A_\infty$-category} $\cC$ 
over a base field $K$ 
consists of 
a set of objects $\Ob(\cC)=\{X,Y,\dots\}$, 
a $\Z$-graded $K$-vector space $\cC(X,Y)=\oplus_{r\in\Z} \cC^r(X,Y)$ for 
each two objects $X,Y\in\Ob(\cC)$ and a collection of 
multilinear maps 
\begin{equation*}
 \m:=\{m_n:\cC(X_1,X_2)\otimes\cdots\otimes \cC(X_n,X_{n+1})
 \to \cC(X_1,X_{n+1})\}_{n\ge 1}
\end{equation*}
of degree $(2-n)$ satisfying the $A_\infty$-relations (\ref{Ainfty-rel}). 

In particular, an $A_\infty$-category $\cC$ with 
higher products all zero, $m_3=m_4=\cdots=0$, is called 
a {\em DG category}. 
\end{defn}
\begin{defn}\label{defn:coh-Ainftycat}
For an $A_\infty$-category $\cC$, a (non-unital) graded category $H(\cC)$ 
is defined by $\Ob(H(\cC)):=\Ob(\cC)$ 
and for any $X,Y\in\cC$ the space of morphisms is 
the cohomology of the chain complex $(\cC(X,Y),m_1)$: 
\begin{equation*}
 H(\cC)(X,Y) :=H(\cC(X,Y),m_1). 
\end{equation*}
The composition in $H(\cC)$ is given 
by the one induced from $m_2$ in $\cC$. 
We call $H(\cC)$ the {\em cohomology of $\cC$}. 

We also define a (non-unital) category $H^0(\cC)$ by 
$\Ob(H^0(\cC)):=\Ob(\cC)$ 
and for any $X,Y\in\cC$ the space of morphisms is 
the cohomology of degree zero only: 
\begin{equation*}
 H^0(\cC)(X,Y):=H^0(\cC(X,Y),m_1). 
\end{equation*}
The composition in $H^0(\cC)$ is given 
by the one induced from $m_2$ in $\cC$. 
We call $H^0(\cC)$ the {\em zero-th cohomology of $\cC$}. 
\end{defn}
The suspension $s\cC$ of an $A_\infty$-category $\cC$ 
is defined by the shift 
\begin{equation*}
 s: \cC(X,Y)\to s(\cC(X,Y))=:(s\cC)(X,Y)
\end{equation*}
for any $X,Y\in\Ob(\cC)=\Ob(s\cC)$, 
where the degree $|b_n|$ of the $A_\infty$-products becomes one 
for all $n\ge 1$ as it does for $A_\infty$-algebras. 
\begin{defn}[$A_\infty$-functor]
Given two $A_\infty$-categories $\cC$, $\cD$, 
an {\em $A_\infty$-functor} 
$\f:=\{f;f_1,f_2,\dots\}:\cC\to \cD$ is 
a map $f:\Ob(s(\cC))\to\Ob(s(\cD))$ of objects with 
degree preserving multilinear maps 
\begin{equation*}
 f_k: (s\cC)(X_1,X_2)\otimes\cdots\otimes (s\cC)(X_k,X_{k+1})
 \to (s\cD)(f(X_1),f(X_{k+1}))
\end{equation*}
for $k\ge 1$ satisfying the defining relations 
of an $A_\infty$-morphism (\ref{Ainfty-morp}). 

In particular, if $f:\Ob(s\cC)\to\Ob(s\cD)$ is bijective 
and $f_1: (s\cC)(X,Y)\to (s\cD)(f(X),f(Y))$ 
induces an isomorphism between 
the cohomologies for any $X,Y\in\Ob(s\cC)$, 
we call the $A_\infty$-functor $\f$
an {\em $A_\infty$-quasi-isomorphism functor}. 
If $f_1$ is itself an isomorphism, then $\f$ is called 
an {\em $A_\infty$-isomorphism functor}. 
\end{defn}

 \subsection{Homological perturbation theory (HPT)}
\label{ssec:hpt}

For an $A_\infty$-algebra $A=(A,\m^A)$,
$(A,d^A:=m^A_1)$ forms a complex. 
The HPT starts with what is called {\em strong deformation retract (SDR)} data 
\cite{gugjds,GLS:chen,GLS,HK}: 
\begin{equation}\label{SDR}
\Nsddata{B}{\ \iota}{\ \pi}{A}{h}, 
\end{equation}
where $(B, d^B:=\pi\circ d^A\circ\iota)$ is a complex with 
chain maps $\iota$ and $\pi$ 
so that $\pi\circ\iota=\Id_B$ and 
$h:A\to A$ is a contracting homotopy satisfying 
\begin{equation}\label{SDR}
 d^Ah+hd^A=\Id_A-P, \qquad P:=\iota\circ\pi. 
\end{equation}
By this definition, $P$ is an idempotent in $A$ 
which commutes with $d^A$, $P d^A=d^A P$. 
If $d^A P=0$, then the SDR (\ref{SDR}) gives a Hodge 
decomposition of the complex $(A,d^A)$, where 
$P(A)\simeq H(A)$ gives the cohomology. 

For an $A_\infty$-algebra $(A[1],\b^A)$ and an SDR data
$\Nsddata{B}{\ \iota}{\ \pi}{A}{h}$, 
there exists an explicit construction of 
an $A_\infty$-algebra $(B[1],\b^B)$ and an $A_\infty$-quasi-isomorphism 
$\f:(B[1],\b^B)\to (A[1],\b^A)$. 
First we set $f_1:=\iota: B[1]\to A[1]$. 
The $A_\infty$-quasi-isomorphism $\f$ is then defined 
inductively by 
\begin{equation*}
  f_{n+1}:= -h\sum_{j\ge 2}\sum_{k_1+\cdots + k_j=n+1}
  b^A_j (f_{k_1}\otimes\cdots\otimes f_{k_j})
\end{equation*}
for given $f_1,f_2,\dots, f_n$. 
Then, the $A_\infty$-structure $\b^B$ in $B[1]$ is given by 
$b^B_1=d^B$ and the formula 
\begin{equation}\label{hpt}
 b^B_{n+1}:= 
 \pi\sum_{j\ge 2}\sum_{k_1+\cdots + k_j=n+1}
 b^A_j (f_{k_1}\otimes\cdots\otimes f_{k_j})
\end{equation}
for $n=1,2,\dots$. 
\begin{lem}\label{lem:hpt}
This $(B[1],\b^B)$ actually forms an $A_\infty$-algebra
and this $\f:(B[1],\b^B)\to (A[1],\b^A)$ is an $A_\infty$-quasi-isomorphism. 
\qed\end{lem}
This formula is given in \cite{GLS,HK} when $(A,\m^A)$ is a DG algebra, 
and in \cite{KoSo}, for instance, for the general case.
The statement of Lemma \ref{lem:hpt} is 
also called the homotopy transfer lemma, see \cite{loday-vallette}, 
which is actually a part of the HPT developed for instance 
by \cite{gugjds,GLS:chen,GLS,HK}. 

As a corollary of Lemma \ref{lem:hpt}, 
when we start from an SDR data defining the Hodge decomposition of
the complex $(A,d^A)$ where $B\simeq H(A)$, 
one obtains the minimal model theorem by
Kadeishvili \cite{kadei:minimal}: 
{\em 
for a given $A_\infty$-algebra $(A,\m^A)$, 
there exists a minimal $A_\infty$-algebra structure $\m^{H(A)}$ in
$H(A)$ and 
an $A_\infty$-quasi-isomorphism $\f:(H(A),\m^{H(A)})\to (A,\m)$. }
Here, that $\m^{H(A)}$ is minimal means $m^{H(A)}_1=0$. 
Such a minimal $A_\infty$-algebra $(H(A),\m^{H(A)})$ is called 
a {\em minimal model} of $(A,\m^A)$. 
We see that the cohomology algebra $H(A,\m^A)$ is obtained 
from a minimal model $(H(A),\m^{H(A)})$ of $(A,\m)$
by forgetting the higher $A_\infty$-products $m^{H(A)}_3,m^{H(A)}_4,\dots$. 

It is straightforward to generalize this HPT for $A_\infty$-algebras to
that for $A_\infty$-category $\cC$ \cite{KoSo}.
We begin with an SDR 
\begin{equation*}
 \Nsddata{\cD(X,Y)}{\ \iota_{XY}}{\ \pi_{XY}}{\cC(X,Y)}{h_{XY}}, 
\end{equation*}
for any $X,Y\in\cC$. 
Then the $A_\infty$-structure of $\cD$ is given by
the straightforward generalization of the formula (\ref{hpt}) 
and an $A_\infty$-quasi-isomorphism functor
\begin{equation*}
  \f:\cD\to\cC
\end{equation*}  
is also given as above.

 \section{The shifts and the signs in $A_\infty$-categories}
\label{sec:main}

Before defining $\ti{\cC}^{(1)}$ and $\ti{\cC}^{(2)}$ mentioned in the
introduction, 
in subsection \ref{ssec:shift}, 
we begin with constructing $A_\infty$-categories
$\cC'=\cC^{(a)}$, $a=1,2$, from an $A_\infty$-category $\cC$. 
We extend this to $\ti{\cC}^{(a)}$, $a=1,2$, in subsection \ref{ssec:main},
and obtain 
Proposition \ref{prop:morp-preserve} and Theorem \ref{thm:HPT-preserve} 
by just 
rewriting Proposition \ref{prop:morp-preserve'} and Theorem \ref{thm:HPT-preserve'} obtained in subsection \ref{ssec:shift}.

 \subsection{The shifts and the signs}
\label{ssec:shift}

For an $A_\infty$-category $\cC$ and an object $X\in\cC$, 
we would like to define another  $A_\infty$-category $\cC'$
where $X\in\cC$ is replaced by $X[1]$ in the sense that
\begin{equation*}
  (\cC')^r(Y,X[1]):=\cC^{r+1}(Y,X) ,\qquad
  (\cC')^r(X[1],Y):=\cC^{r-1}(X,Y) . 
\end{equation*}
We of course set $\cC'(Y,Z):=\cC(Y,Z)$ if $Y\ne X$ and $Z\ne X$. 
For a morphism $\alpha$ in $\cC$, we denote by $\alpha'$
the corresponding morphism in $\cC'$. 

A way of defining an $A_\infty$-structure in $\cC'$
is presented explicitly in \cite{fukaya:mirror2}.
For 
\begin{equation*}
 \begin{array}{cccc}
   b_k:& s\cC(X_{1},X_{2})\times\cdots\times s\cC(X_{k},X_{k+1})
   & \to & s\cC(X_{1},X_{k+1}) \\
   & (\bar{\alpha}_{12},\dots,\bar{\alpha}_{k(k+1)})& \mapsto &
   b_k(\bar{\alpha}_{12},\dots,\bar{\alpha}_{k(k+1)}) 
 \end{array}
\end{equation*}
in $\cC$, we set $b'_k$ to be
\begin{equation*}
  b'_k(\bar{\alpha}'_{12},\dots,\bar{\alpha}'_{k(k+1)})= \pm
  ( b_k(\bar{\alpha}_{12},\dots,\bar{\alpha}_{k(k+1)}) )', 
\end{equation*}
where the sign $\pm$ is $-1$ when $X_{1}=X$ and $+1$ otherwise. 
One can check directly that this $\m'$ preserves 
the $A_\infty$-relations of $\m$ \cite{fukaya:mirror2}. 
We denote the resulting $A_\infty$-category by $\cC'=\cC^{(1)}$. 

Though the definition of this $\m'=:\m^{(1)}$ is natural enough, 
there are other choices of $A_\infty$-structures. 
One such $A_\infty$-structure is given by 
\begin{equation*}
  b'_k(\bar{\alpha}'_{12},\dots,\bar{\alpha}'_{k(k+1)})= \pm
  ( b_k(\bar{\alpha}_{12},\dots,\bar{\alpha}_{k(k+1)}) )', 
\end{equation*}
where the sign $\pm$ is $-1$ if the number of $X$ in $X_{1},\dots, X_{k}$ 
(, not $X_{1},\dots, X_{k+1}$,) 
is odd and $+1$ if it is even. 
We denote the resulting $A_\infty$-category by $\cC^{(2)}$. 
In general, these two $A_\infty$-categories are not 
linearly $A_\infty$-isomorphic to each other. 

First, by direct calculations, we see that the following holds. 
\begin{prop}\label{prop:morp-preserve'}
Any $A_\infty$-functor $\f:\cD\to\cC$ induces 
an $A_\infty$-functor $\f^{(a)}:\cD^{(a)}\to\cC^{(a)}$, with $a=1$ or $2$,
given by 
\begin{equation*}
  f_k^{(a)}(\bar{\alpha}'_{12},\dots, \bar{\alpha}'_{k(k+1)})
  := ( f_k(\bar{\alpha}_{12},\dots, \bar{\alpha}_{k(k+1)}) )'.  
\end{equation*}
\end{prop}
\begin{pf}
We can check directly that this satisfies the defining equations of 
$A_\infty$-morphisms (\ref{Ainfty-morp}). 
\qed\end{pf}

Now, we discuss both $\cC^{(1)}$ and $\cC^{(2)}$ are natural enough
by showing that the $A_\infty$-structures of both types are preserved by
the HPT. 
First, assume that an SDR data for $\cC$ is given.
Namely, each space $\cC(Y,Z)$ of morphisms has 
an SDR 
\begin{equation*}
 \Nsddata{\cD(Y,Z)}{\ \iota_{YZ}}{\ \pi_{YZ}}{\cC(Y,Z)}{h_{YZ}}, 
\end{equation*}
so that 
\begin{equation}\label{SDR-shift}
 d_{YZ}h_{YZ}+h_{YZ}d_{YZ}= \Id - \iota_{YZ}\circ\pi_{YZ} . 
\end{equation}
Then, an SDR data is induced for $\cC^{(a)}$ as follows.  
In both cases $a=1$ and $2$,
the $b_1'=b_1^{(a)}$ above satisfies
\begin{equation}\label{SDRinduce1}
  b_1^{(a)}(\bar{\alpha}'):= 
  \begin{cases}
   - (b_1(\bar{\alpha}))' & Y=X \\
   (b_1(\bar{\alpha}))' & \text{otherwise}
  \end{cases}
\end{equation}
for $\bar{\alpha}\in s\cC(Y, Z)$. 
Let us denote by $Y'$ and $Z'$ the objects in $\cC^{(a)}$ 
corresponding to $Y$ and $Z$ in $\cC$.
Namely, 
$Y'=X[1]$ if $Y=X$ and $Y'=Y$ otherwise, and 
$Z'=X[1]$ if $Z=X$ and $Z'=Z$ otherwise. 
Then, (\ref{SDRinduce1}) 
is equivalent to that 
$d_{Y'Z'}^{(a)}=m_1^{(a)}$ satisfies 
\begin{equation*}
  d_{Y'Z'}^{(a)}(\alpha'):= 
  \begin{cases}
   - (d_{YZ}(\alpha))' & Y=X \\
   (d_{YZ}(\alpha))' & \text{otherwise} . 
  \end{cases}
\end{equation*}
Therefore, setting 
\begin{equation*}
  h^{(a)}_{Y'Z'}(\alpha'):= 
  \begin{cases}
   -(h_{YZ}(\alpha))' & Y=X \\
   (h_{YZ}(\alpha))' & \text{otherwise}
  \end{cases} , 
\end{equation*}
the formula (\ref{SDR-shift}) is preserved, 
\begin{equation*}
  d^{(a)}_{Y'Z'}h^{(a)}_{Y'Z'}+h^{(a)}_{Y'Z'}d^{(a)}_{Y'Z'}
  = \Id - \iota_{Y'Z'}\circ\pi_{Y'Z'} ,  
\end{equation*}
where
$\iota_{Y'Z'}$ is defined as $\iota_{Y'Z'}(\alpha'):=(\iota_{YZ}(\alpha))'$
and $\pi_{Y'Z'}$ is defined similarly. 
We call this
the {\em induced SDR data} on $\cC^{(a)}$. 

\begin{thm}\label{thm:HPT-preserve'}
For a given $A_\infty$-category $\cC$, 
we fix an SDR data on $\cC$.
By the HPT, we obtain an $A_\infty$-category $\cD$
and an $A_\infty$-quasi-isomorphism functor 
$\f:\cD\to\cC$. 
Similarly, we obtain an $A_\infty$-quasi-isomorphism functor
$\g:\cE\to\cC^{(a)}$ with $a=1$ or $2$ 
by applying the HPT for $\cC^{(a)}$ with the induced SDR data. 
Then, one has $\cE=\cD^{(a)}$ and $\g=\f^{(a)}$. 
\end{thm}
\begin{pf}
We prove this when $a=1$. 

We first show that $\g=\f^{(1)}$. 
Recall that $\f$ is given inductively by 
$f_1=\iota$ and 
\begin{equation*}
f_i(\cdots)=
-h\sum_{k\ge 2}\sum_{i_1+\cdots +i_k=i}
b^{\cC}_k(f_{i_1}(\cdots),\dots,f_{i_k}(\cdots))
\end{equation*}
for $i\ge 2$. 
The $A_\infty$-quasi-isomorphism $\g$ is defined
in the same way but with replacing $b^\cC_k$ by $b_k^{\cC^{(a)}}$. 
{}For 
\begin{equation*}
\begin{array}{cccc}
  f_k & :s\cD(X_{1},X_{2})\times\cdots\times s\cD(X_{k},X_{k+1})
 &\to & \cC(X_{1},X_{k+1}) \\ 
& (\bar{\alpha}_{12},\dots,\bar{\alpha}_{k(k+1)}) & \mapsto & 
f_k(\bar{\alpha}_{12},\dots,\bar{\alpha}_{k(k+1)})  , 
\end{array}
\end{equation*}
we show by induction that 
\begin{equation}\label{induction1}
g_k(\bar{\alpha}'_{12},\dots,\bar{\alpha}'_{k(k+1)})=
 ( f_k(\bar{\alpha}_{12},\dots,\bar{\alpha}_{k(k+1)}) )'
\end{equation}
holds.

Assuming this holds for $k\le n$,
one has
\begin{equation*}
   g_{n+1}(\cdots)= -h^{(a)}\sum_{k\ge 2}\sum_{i_1+\cdots +i_k=n+1}
b^{\cC^{(a)}}_k(f_{i_1}(\cdots),\dots,f_{i_k}(\cdots)), 
\end{equation*}
where the sign coming from $h^{(a)}$ and that from $b^{(a)}_k$
cancel with each other. 
Namely, (\ref{induction1}) holds true for $k=n+1$. 
Thus, the induction is completed and we have $\g=\f^{(1)}$. 

Next, we show that $b^\cE_k=b^{\cD^{(1)}}_k$. 
We assume this holds for $1\le k\le n$. 
For $k=n+1$, 
\begin{equation*}
\begin{array}{cccc}
  b^\cE_{n+1} & :s\cE(X'_{1},X'_{2})\times\cdots\times s\cE(X'_{n+1},X'_{n+2})
 &\to & s\cE(X'_{1},X'_{n+2}) \\ 
& (\bar{\alpha}'_{12},\dots,\bar{\alpha}'_{(n+1)(n+2)}) & \mapsto & 
b^\cE_{n+1}(\bar{\alpha}'_{12},\dots,\bar{\alpha}'_{(n+1)(n+2)})  , 
\end{array}
\end{equation*}
is given by
\begin{equation*}
  b^\cE_{n+1}(\cdots)=\sum_{k\ge 2}\sum_{i_1+\cdots +i_k=n+1}
  \pi b^{\cC_{(a)}}_k(g_{i_1}(\cdots),\dots, g_{i_k}(\cdots)). 
\end{equation*}
Thus, the minus sign appears from
$b^{\cC_{(1)}}_k(g_{i_1}(\cdots),\dots, g_{i_k}(\cdots))$
if and only if $X_1=X$, i.e., $X'_{1}=X[1]$. 
This shows that $b^\cE_{n+1}=b^{\cD^{(1)}}_{n+1}$. 

The case for $a=2$ can also be shown in the parallel strategy. 
\qed\end{pf}

\subsection{Additive $A_\infty$-categories}
\label{ssec:main}

For a given $A_\infty$-category $\cC$ and an object $X\in\cC$,
by replacing $X$ by $X[1]$ we constructed
the $A_\infty$-category $\cC^{(a)}$ with $a=1$ or $2$. 
We can instead construct an $A_\infty$-category by adding
$X[1]$ or $X[-1]$ to $\cC$.

Repeating this procedure 
yields an additive $A_\infty$-category $\ti{\cC}^{(a)}$ as follows.

In both cases $a=1,2$, 
an object in $\ti{\cC}^{(a)}$ is a finite direct sum 
$\cX:=X_1[r_1]\oplus\cdots\oplus X_l[r_l]$, 
where $X_i\in\cC$ for each $i=1,\dots,l$ and 
$[r_i]$ indicates the formal degree shift by $r_i\in\Z$. 
For $\cX, \cY\in\ti{\cC}$, 
the space $\ti{\cC}(\cX,\cY)=\ti{\cC}^{(a)}(\cX,\cY)$
of morphisms is given by 
\begin{equation*}
 \ti{\cC}^r(X[n],Y[m]):=\cC^{r+m-n}(X,Y) . 
\end{equation*}
The $A_\infty$-structure $\ti{\b}^{(a)}=\ti{\b}^{\ti{\cC}^{(a)}}$
is determined by defining
\begin{equation*}
  \ti{b}_k^{(a)}(\bar{\alpha}'_{12},\dots,\bar{\alpha}'_{k(k+1)})
  \in s\ti{\cC}(X_1[r_1],X_{k+1}[r_{k+1}])
\end{equation*}
for $\bar{\alpha}_{i(i+1)}\in s\ti{\cC}(X_i,X_{i+1})$ and 
$\bar{\alpha}'_{i(i+1)}\in s\ti{\cC}(X_i[r_i],X_{i+1}[r_{i+1}])$. 
Repeating the construction of $\cC^{(a)}$ in the previous subsection, 
we set 
\begin{equation}\label{tib1}
  \ti{b}_k^{(1)}(\bar{\alpha}'_{12},\dots,\bar{\alpha}'_{k(k+1)})
  = (-1)^{r_1} ( b_k(\bar{\alpha}_{12},\dots,\bar{\alpha}_{k(k+1)}) )'
\end{equation}
and 
\begin{equation}\label{tib2}
  \ti{b}_k^{(2)}(\bar{\alpha}'_{12},\dots,\bar{\alpha}'_{k(k+1)})
  = (-1)^{r_1+\cdots + r_k}
  ( b_k(\bar{\alpha}_{12},\dots,\bar{\alpha}_{k(k+1)}) )'. 
\end{equation}
Then, both $\ti{\cC}^{(1)}$ and $\ti{\cC}^{(2)}$ form $A_\infty$-categories. 

In the unsuspended notation, by (\ref{m-suspension}) 
we have 
\begin{equation}\label{tim1}
  \ti{m}_k^{(1)}(\alpha_{12}',\dots,\alpha_{k(k+1)}')
  = (-1)^{r_1k+r_2+\cdots + r_k} ( m_k(\alpha_{12},\dots,\alpha_{k(k+1)}) )'
\end{equation}
and 
\begin{equation}\label{tim2}
  \ti{m}_k^{(2)}(\alpha_{12}',\dots,\alpha_{k(k+1)}')
  = (-1)^{r_1k} ( m_k(\alpha_{12},\dots,\alpha_{k(k+1)}) )',  
\end{equation}
where $\alpha_{i(i+1)}\in\cC(X_i,X_{i+1})$. 

As corollaries of the previous subsection, we obtain the followings. 
We let $a$ be $1$ or $2$. 
\begin{prop}\label{prop:morp-preserve}
An $A_\infty$-functor $\f:\cD\to\cC$ induces
an $A_\infty$-functor $\f^{(a)}:\ti{\cD}^{(a)}\to\ti{\cC}^{(a)}$
which is defined by 
\begin{equation*}
  f_k^{(a)}(\bar{\alpha}'_{12},\dots,\bar{\alpha}'_{k(k+1)})
 := ( f_k(\bar{\alpha}_{12},\dots,\bar{\alpha}_{k(k+1)}) )' . 
\end{equation*}
\qed\end{prop}
\begin{thm}\label{thm:HPT-preserve}
For a given $A_\infty$-category $\cC$, 
we fix an SDR data on $\cC$. 
By HPT, we obtain an $A_\infty$-category $\cD$
and an $A_\infty$-quasi-isomorphism functor 
$\f:\cD\to\cC$, which further induces an $A_\infty$-quasi-isomorphism
functor $\ti{\f}^{(a)}:\ti{\cD}^{(a)}\to\ti{\cC}^{(a)}$. 
Similarly, we obtain an $A_\infty$-quasi-isomorphism functor
$\ti{\g}:\ti{\cE}\to\ti{\cC}^{(a)}$ 
by applying HPT for $\ti{\cC}^{(a)}$ with the induced SDR data.
Then one has 
\begin{equation*}
 \ti{\cE}=\ti{\cD}^{(a)} ,\qquad  \ti{\g}=\ti{\f}^{(a)} . 
\end{equation*}
\qed\end{thm}
This theorem indicates 
that the procedure of taking $\ti{\cC}^{(a)}$ from $\cC$
is compatible with the HPT.

 \subsection{Triangulated $A_\infty$-categories}
\label{ssec:tri}

Once we construct $\ti{\cC}^{(a)}$,
we can further construct an $A_\infty$-category
$\Tw^{(a)}(\cC)$ of one-sided twisted complexes. 
When $\cC$ is strictly unital, then
$\Tw^{(a)}(\cC)$ is also strictly unital, and 
its zero-th cohomology 
\begin{equation*}
  \Tri^{(a)}(\cC):=H^0(\Tw^{(a)}(\cC))  
\end{equation*}
forms a triangulated category 
as proposed in \cite{HMS}. 
In this sense, $\Tw^{(a)}(\cC)$ is called a triangulated
$A_\infty$-category. 
The construction of $\Tw^{(a)}(\cC)$ is discussed explicitly
in \cite{fukaya:mirror2}. 
When the shift functor $T:\Tri^{(a)}(\cC)\to \Tri^{(a)}(\cC)$ is given,
the exact triangles are defined quite naturally
as we explain later. 

In this subsection, 
we present the shift functor $T:\Tri^{(a)}(\cC)\to \Tri^{(a)}(\cC)$
explicitly. 

First, we explain the triangulated $A_\infty$-category $\Tw^{(a)}(\cC)$ 
constructed explicitly in \cite{fukaya:mirror2}. 
For the notations, we mainly follow \cite{hk:triAinfty}. 

A {\em one-sided twisted complex} $(\cX,\Phi)$ is
a pair of an object $\cX:=X_1[r_1]\oplus\cdots\oplus X_l[r_l]\in\ti{\cC}^{(a)}$ and a degree zero morphism 
$\Phi=\oplus_{i,j=1}^l \phi_{ij}\in\ti{\cC}^{(a)}(\cX,\cX)$
which satisfies 
$\phi_{ij}=0$ for $i\ge j$ and  
the $A_\infty$-Maurer-Cartan equation:
\begin{equation}\label{MC}
  \ti{b}_1^{(a)}(\Phi)+\ti{b}_2^{(a)}(\Phi,\Phi) + \cdots =0 . 
\end{equation}
The category $\Tw^{(a)}(\cC)$ consists of one-sided twisted complexes,
where the spaces of morphisms are defined by
$\Tw^{(a)}(\cC)((\cX,\Phi),(\cY,\Psi)):=\ti{\cC}^{(a)}(\cX,\cY)$, and 
the $A_\infty$-structure $\b^{\Tw^{(a)}}$ is given by 
\begin{equation*}
 \begin{split}
 &  b_n^{\Tw^{(a)}}(\varphi_{12},\dots,\varphi_{n(n+1)}) \\
     & =\sum_{k_1,\dots,k_{n+1}\ge 0}
     \ti{b}^{(a)}_{n+k_1+\cdots +k_{n+1}}
     ((\Phi_1)^{k_1},\varphi_{12},(\Phi_2)^{k_2},\dots, \varphi_{n(n+1)},(\Phi_{n+1})^{k_{n+1}})  
 \end{split}
\end{equation*}
for $(\cX_i,\Phi_i)\in\Tw^{(a)}(\cC)$, $i=1,\dots, n+1$. 
Then, $\Tw^{(a)}(\cC)$ again forms an
$A_\infty$-category \cite{fukaya:mirror2}.

Next, we define an additive isomorphism
$T: \ti{\cC}^{(a)}\to\ti{\cC}^{(a)}$ 
which satisfies 
\begin{equation*}
 T(X)=X[1],\qquad X\in\cC\subset\ti{\cC}^{(a)}
\end{equation*}
and 
\begin{equation}\label{shift-functor-Ainfty}
 T\, \ti{b}_k^{(a)} = (-1)^k\ti{b}_k^{(a)}(T\otimes\cdots\otimes T) . 
\end{equation}
\footnote{Unfortunately, this $T$ does not form a linear $A_\infty$-automorphism
$T: \ti{\cC}^{(a)}\to\ti{\cC}^{(a)}$ because of the sign $(-1)^k$ 
in (\ref{shift-functor-Ainfty}). }
{}For $\alpha_{ij}\in\cC(X_i,X_j)$, 
denote the corresponding morphisms by
$\alpha'_{ij}\in\ti{\cC}^{(a)}(X_i[r_i],X_j[r_j])$ and
$\alpha''_{ij}\in\ti{\cC}^{(a)}(X_i[r_i+1],X_j[r_j+1])$. 
Then, we set
$T:\ti{\cC}^{(1)}(X_i[r_i],X_j[r_j])\to\ti{\cC}^{(1)}(X_i[r_i],X_j[r_j])$ 
by 
\begin{equation*}
  T(\alpha'_{ij})=-\alpha''_{ij} , 
\end{equation*}
and $T:\ti{\cC}^{(2)}(X_i[r_i],X_j[r_j])\to\ti{\cC}^{(2)}(X_i[r_i],X_j[r_j])$ 
by
\begin{equation*}
  T(\alpha'_{ij})= \alpha''_{ij} . 
\end{equation*}
These determine $T$ for all objects and morphisms in $\ti{\cC}^{(a)}$,
and we see that (\ref{shift-functor-Ainfty}) actually holds. 

This $T$ naturally induces an additive isomorphism 
$T: \Tw^{(a)}(\cC)\to\Tw^{(a)}(\cC)$ as follows. 
For objects, in both cases $a=1$ and $2$,
we first observe that $-T(\Phi)$ satisfies
the $A_\infty$-Maurer-Cartan equation
if $\Phi$ does 
due to (\ref{shift-functor-Ainfty})
and (\ref{tib1}) or (\ref{tib2}). 
Thus, we set
\begin{equation*}
  T(\cX,\Phi):=(T(\cX), -T(\Phi)), 
\end{equation*}
where the $T$'s in the left hand side are those in $\ti{\cC}^{(a)}$. 
For morphisms, for
$\varphi_{ij}\in\Tw^{(a)}(\cC)((\cX_i,\Phi_i),(\cX_j,\Phi_j))$, 
we set $T(\varphi_{ij})$ just as $T(\varphi_{ij})$ in the sense in $\ti{\cC}^{(a)}$. 
Then, again by (\ref{shift-functor-Ainfty}), 
we have
\begin{equation*}
  T\, b_k^{\Tw^{(a)}(\cC)} = (-1)^k b_k^{\Tw^{(a)}(\cC)}( T\otimes\cdots\otimes T) . 
\end{equation*}
In particular, $T m_2^{\Tw^{(a)}(\cC)} =  m_2^{\Tw^{(a)}(\cC)}( T\otimes T)$ 
holds for $k=2$, so this induces
an additive isomorphism functor $T:\Tri^{(a)}(\cC)\to\Tri^{(a)}(\cC)$.  

We treat this $T$ as the shift functor.
For an $m_1^{\Tw^{(a)}(\cC)}$-closed morphism 
$\varphi\in \Tw^{(a)}(\cC)((\cX,\Phi_\cX), (\cY,\Phi_\cY))$
of degree zero, 
the mapping cone $C(\varphi)$ is defined by 
\begin{equation*}
  C(\varphi):=\left(T(\cX)\oplus\cY,
  \left( \bps -T(\Phi_\cX)& \varphi' \\ 0 & \Phi_\cY \eps\right)\right) , 
\end{equation*}
where $\varphi'$ is the morphism in
$s\Tw^{(a)}(\cC)(T(\cX,\Phi_\cX), (\cY,\Phi_\cY))$
induced from $\varphi$. 
Then, 
exact triangles in $\Tri^{(a)}(\cC)$ are defined as sequences
which are isomorphic to 
\begin{equation*}
  \cdots (\cX,\Phi_\cX)\overset{\varphi}\to (\cY,\Phi_\cY)\to
  C(\varphi)\to T(\cX,\Phi_\cX)\to\cdots 
\end{equation*}
with some $\varphi$. 
As mentioned in \cite{fukaya:mirror2}, 
it is shown by direct calculations that
these actually satisfy the axiom of exact triangles. 
Thus, $\Tri^{(a)}(\cC)$ forms a triangulated category
when $\cC$ is a strictly unital $A_\infty$-category.

 \subsection{Concluding remarks}
\label{ssec:consequence}

For a given $A_\infty$-functor $\f:\cD\to\cC$ of strictly unital
$A_\infty$-categories, 
the $A_\infty$-functor $\ti{\f}^{(a)}:\ti{\cD}^{(a)}\to\ti{\cC}^{(a)}$ 
in Proposition \ref{prop:morp-preserve} 
further induces
an $A_\infty$-functor $\Tw^{(a)}(\cD)\to\Tw^{(a)}(\cC)$
and then a functor $\Tri^{(a)}(\cD)\to\Tri^{(a)}(\cC)$
of triangulated categories 
as explained in \cite{fukaya:mirror2, Seidel:book}, etc. 
In particular, if $\f$ is an $A_\infty$-equivalence functor, then
the induced functor between the triangulated categories
is an equivalence. See \cite{Seidel:book}.

Note that $\ti{\cC}_1$ and $\ti{\cC}_2$ are not $A_\infty$-isomorphic
to each other in general.
Actually, there exist such examples of $\cC$ even if $m_1=m_3=\cdots =0$. 
This means, 
there does not exist any $A_\infty$-isomorphism functor
between 
$\Tw^{(1)}(\cC)$ and $\Tw^{(2)}(\cC)$
so that it preserves $\ti{\cC}^{(a)}\subset\Tw^{(a)}(\cC)$. 
In this sense, 
there does not exist any natural isomorphism functor of triangulated categories
between 
$\Tri^{(1)}(\cC)$ and $\Tri^{(2)}(\cC)$.

%\bibliographystyle{plain}
%\bibliography{cite}
%\bibliography{cite-ag}

\end{document}